\def\marker{\>\hbox{${\vcenter{\vbox{
    \hrule height 0.4pt\hbox{\vrule width 0.4pt height 6pt
    \kern6pt\vrule width 0.4pt}\hrule height 0.4pt}}}$}\>}

\documentclass[12pt]{article}
\usepackage{fullpage}
\usepackage{amsmath,amsfonts,amssymb,amsthm,graphics,amscd,MnSymbol}
\usepackage{graphicx}
\usepackage{subcaption}

\usepackage[
pdfauthor={},
pdftitle={ Strong list-chromatic index of subcubic graphs is at most $10$},
pdfstartview=XYZ,
bookmarks=true,
colorlinks=true,
linkcolor=blue,
urlcolor=blue,
citecolor=blue,
bookmarks=false,
linktocpage=true,
hyperindex=true
]{hyperref}

\newtheorem {Theorem}  {Theorem}[section]
\newtheorem {Lemma}[Theorem]{Lemma}

\newtheorem {Corollary}[Theorem]{Corollary}
\newtheorem {Claim}{Claim}
\theoremstyle{definition}
\newtheorem{Definition}[Theorem]{Definition}

\newenvironment {Proof} {\noindent {\bf Proof.}}{\quad $\square$\par\vspace{3mm}}

\date{\today}
\title{\bf Strong list-chromatic index of subcubic graphs is at most $10$}

\author{Yunfang Tang\thanks{
Department of Mathematics, China Jiliang University, Hangzhou 310018, P.R. China. Grant number:
  NSFC11701543. Email: tangyunfang8530@cjlu.edu.cn.}
\and Zhiwei Bi \thanks{
Department of Mathematics, China Jiliang University, Hangzhou 310018, P.R. China. Email: S23080701004@cjlu.edu.cn.}
\\[0.2cm]
      }

\begin{document}

\maketitle

\vspace{-2pc}

\begin{abstract}
A strong edge coloring of a graph $G$ is an assignment of colors to the edges of $G$ such that two distinct edges are colored differently if they are incident to a common edge or share an endpoint. The strong chromatic index of a graph $G$, denoted by $\chi_{s}'(G)$, is the minimum number of colors needed for a strong edge coloring of $G$.  The edge weight of a graph $G$ is defined to be $\max\limits_{uv\in E(G)}\{(d_G(u)+d_G(v))\}$. It was proved in Chen et al in 2020 that every graph with edge weight at most 6 has a strong edge-coloring using at most 10 colors. In this paper, we consider the list version of strong edge-coloring. We strengthen this result by showing that every graph with edge weight at most 6 has a strong list-chromatic index at most 10. Specially, every subcubic graph has a strong list-chromatic index at most 10, which improves a result of Dai et al. in 2018.

\end{abstract}

\section{\textbf{Introduction}}

All graphs in this paper are finite and simple. Let $G$ be a graph with vertex set $V(G)$ and edge set $E(G)$. For positive integer $k$, let $[k]=\{1,2,\ldots,k\}$. If there is a mapping $\phi:E(G) \rightarrow [k]$ such that any two edges meeting at a common vertex, or joining by the same edge of $G$, are assigned different values (colors), then we call $\phi$ a \textit{strong $k$-edge coloring}, and $G$ is \textit{strong $k$-edge colorable}. The \textit{strong chromatic index} of $G$, denoted by $\chi_{s}'(G)$, is the smallest number $k$ for which $G$ has a strong $k$-edge coloring. For each edge $e$ of $G$, let $L(e)$ be the set of all colors for $e$, and let $L$ be a \textit{k-edge list assignment} of $G$ if $|L(e)|\ge k$ for all $e\in E(G)$. For any $k$-edge list assignment $L$ of $G$, if there exists a strong edge coloring $\phi$ such that $\phi(e)\in L(e)$ for all $e\in E(G)$, then $G$ is \textit{strongly $k$-edge choosable} and $\phi$ is called a \textit{strong $L$-edge-coloring}. The \textit{strong list-chromatic index}, denoted by $\chi_{s,l}'(G)$, is the minimum positive integer $k$ for which $G$ is strongly $k$-edge choosable. Note that $\chi_s'(G)\leq\chi_{s,l}'(G)$ for every graph $G$.

For a graph $G$, for $u\in V(G)$, let $d_G(u)$ denote the degree of $u$, and an \textit{ $i$-cycle} is a cycle of length $i$ in $G$, denoted by $C_n$. Let $C_n^+$ be a graph obtained from $C_n$ by joining a pendent vetex to each vertex of $C_n$, respectively.
The\textit{ edge weight} of a graph $G$ is defined to be $max\{d_G(u)+d_G(v):uv\in E(G)\}$. Given two edges $e,e'\in E(G)$, we say that \textit{$e$ sees $e'$} if $e$ and $e'$ are adjacent, or share a common adjacent edge, i.e., $d_G(e,e')\le 2$. So, equivalently, for any $k$-edge list assignment $L$, if there exists an assignment $\phi$ with $\phi(e)\in L(e)$ for all $e\in e(G)$ such that every two edges that can see each other receive distinct colors, then $G$ is strongly $k$-edge choosable.

In 2018, Dai et al~\cite{Dai2018} proved  that if $G$ a subcubic graph, $\chi_{s,l}'(G) \leq 11$. In 2020, Chen et al~\cite{Chen2020} proved that if $G$ is a graph with the edge weight at most 6, then $\chi_s'(G) \leq 10$. This paper improves and extends those results by proving that if $G$ is a subcubic graph, then $\chi_{s,l}'(G) \leq 10$ (see Theorem~\ref{key}); if a graph $G$ whose edge weight is at most $6$, then $\chi_{s,l}'(G) \leq 10$ (see Corollary~\ref{weight}). The proof in~\cite{Dai2018} mainly uses Hall theorem and the Combinatorial Nullstellensatz.

\section{Key idea}
In this paper, we also use algebraic method Combinatorial Nullstellensatz to find more graph structure properties.

\begin{Definition}\label{TY} Let $\mathbb{F}$ be an arbitrary field and  $J = x_{j_{1}}^{i_1}\cdots x_{j_{k}}^{i_{k}}$. Define a mapping $\eta_J:\mathbb{F}[x_1,x_2,\ldots,x_n]\rightarrow\mathbb{F}[\{x_1,x_2,\ldots,x_n\}\setminus\{x_{j_1},x_{j_2},\ldots,x_{j_k}\}]$ as follows,
for every $P\in \mathbb{F}[x_1,x_2,\ldots,x_n]$, \textit{the coefficient} of the monomial $J$ in $P$,
$$\eta_J[P]=\frac{1}{i_1!\cdots i_k!}\cdot\frac{\partial P}{\partial x_{j_1}^{i_1}\cdots x_{j_k}^{i_k}}\mid_{x_{j_1}=\cdots = x_{j_k}=0}.$$
\end{Definition}

\begin{Theorem}$\mathrm{(}$ Combinatorial Nullstellensatz {\rm\cite{Alon1999}}$\mathrm{)}$ \label{CN} Let $\mathbb{F}$ be an arbitrary field, and let $P = P(x_1, x_2, \ldots, x_n)$ be a polynomial in $\mathbb{F}[x_1, x_2, \ldots, x_n]$. Suppose that the degree $\deg(P)$ of $P$ equals $\sum_{i=1}^n k_i$, where each $k_i$ is a non-negative integer, and the coefficient of $x_1^{k_1} x_2^{k_2} \cdots x_n^{k_n}$ in $P$ is non-zero. Then if $S_1, S_2, \ldots, S_n$ are subsets of $\mathbb{F}$ with $|S_i| > k_i$, $i = 1, 2, \ldots, n$, there exist $s_1 \in S_1, s_2 \in S_2, \ldots, s_n \in S_n$ so that $P(s_1, s_2, \ldots, s_n) \neq 0$.
\end{Theorem}

Let $L$ be a $k$-edge list assignment of $G$. Assume $H$ is an induced subgraph. Let $G'=G-V(H)$ (not necessary connected). By the minimality of $G$, there exists a strong $L$-edge-coloring $\varphi$ of $G'$. We may assume $|L(e)|=k$ for each $e \in E(G)$. We will extend the partial coloring $\varphi$ of $G'$ to a strong $L$-edge-coloring of the whole graph $G$. For convenience, denote by $E_0=E(G)\setminus E(G')$ and $t=|E_0|$.
For $e\in E_0$, let $F_{\varphi}(e)=\{\varphi(e'):e'\in E(G'), dis_G(e,e')\le 2\}$. Set $S(e)=L(e) \setminus F_{\varphi}(e)$ and $s(e)=|S(e)|=k-|F_{\varphi}(e)|$. Let $t=|E_0|$. For $i\in [t]$,and $e_i\in E_0$, assign a variable $x_{i}$  to the edge $e_i$.
Consider a  polynomial $$P_{H}=P(x_1,x_2,\cdots,x_t)=\prod\limits_{\substack{i<j,\\~d_{G[E_0]}(e_i,e_j)\le 2}}x_i-x_j.$$
By the definition of $P_H$, if there exists a mapping $\phi$ with $\phi\mid_{E(G')}=\varphi$, and $\phi(e_i)=s_i$, where  $s_i \in S(e_i)$ for $e_i \in E_0 $ such that $P_{H}(\phi)=P(s_1, s_2,\cdots,s_t) \neq 0 $, then $G$ is strongly $k$-edge choosable, a contradiction. Thus $G$ contains no such subgraph $H$ in $G$. Moreover, by Definition~\ref{TY} and Theorem~\ref{CN}, we have the following corollary.

\begin{Corollary}\label{CN1}
If $P_H$ has a monomial $J=\prod_{e\in E_0}x_e^{i_e}$ such that $\eta_{J}[P_H]\neq 0$, where $\sum \limits_{e\in E_0}i_e =\deg(P_H)$, then $G$ is strongly $k$-edge choosable, where $1+max\{{i_e}:e\in E_0\}\le s(e)$.
\end{Corollary}

In the proof of Theorem~\ref{key}, $G$ is a cubic graph and  $H=C_n=v_1v_2\cdots v_nv_1$ with $n\ge 6$ in $G$ and $k=10$. Note that $E_0=E(C_n^+)$. For convenience, denoted by $E_G(v_i)=\{e_{i-1},e_i,f_i\}$, ($e_0=e_n$), where $E(C_n)=\{e_i|i\in[n]\}$. Assign variables $x_{i}$ and $y_{i}$ to $e_i$ and $f_i$, respectively, for $i\in[n]$. Forthermore, all the coefficients of the monomials we need are specific and can be computed by MATHEMATICAL.

\section{Known results}
This section lists some known results that are needed in the proof of Theorem~\ref{key} or Corollary~\ref{weight}.

The following lemma is proved in~\cite{Dai2018}, which will be used in Theorem~\ref{key}.

\begin{Lemma}{\rm\cite{Dai2018}}\label{useful}
Let $G$ be a subcubic graph, and $L$ be a $10$-edge list assignment of $G$.
If $G$ has no $L$-edge-coloring but any proper subgraph of $G$ has a $L$-edge-coloring,
then $G$ is cubic and has no cycles of length at most 5.
\end{Lemma}

\begin{Theorem}\rm{\cite{Chen2020}}\label{Chen}
Let $G$ be a graph with edge weight at most $6$. Then $\chi_s'(G) \leq 10$.
\end{Theorem}

The following Lemma is in the proof of Theorem~\ref{Chen}.
\begin{Lemma}\rm{\cite{Chen2020}}\label{useful1}
Let $G$ be a graph with edge weight at most $6$, and $L$ be a $10$-edge list assignment of $G$.
If $G$ has no $L$-edge-coloring but any proper subgraph of $G$ has a $L$-edge-coloring,
then $G$ is cubic or $G$ is $(2,4)$-bipartite graph.
\end{Lemma}

\begin{Theorem}\rm{\cite{Deng}}\label{useful2}
Let $G$ be a $(2,\Delta)$-bipartite graph where $\Delta(G) \geq 4$. Then $\chi_{s,l}'(G) \leq 3\Delta-3$.
\end{Theorem}

The main result of this paper is the following theorem.
\begin{Theorem} \label{key}
Let $G$ be a subcubic graph. Then $ \chi_{s,l}'(G) \leq 10 $.
\end{Theorem}

Thus Combining Lemma~\ref{useful1} and Theorems~\ref{useful2} and~\ref{key}, we obtain the following
result:
\begin{Corollary} \label{weight}
Let $G$ be a graph with edge weight at most $6$. Then $ \chi_{s,l}'(G) \leq 10 $.
\end{Corollary}

\section{The proof of Theorem \ref{key}}
Let $G$ be a minimal counterexample such that $E(G)$ is as small as possible.
By Lemma~\ref{useful}, we only consider $G$ is cubic and the girth of $G$ is at least $6$.
Suppose that $G$ is a cubic graph, and $L$ is a $10$-edge list assignment of $G$.
In the following, we will show that $G$ has no $n$-cycles with $n\ge 6$.

\begin{Claim}
$G$ has no $6$-cycles, i.e. $n\ne 6$.
\end{Claim}

\begin{figure}[ht]
    \centering
    \begin{subfigure}[b]{0.3
\textwidth}
        \includegraphics[width=\linewidth]{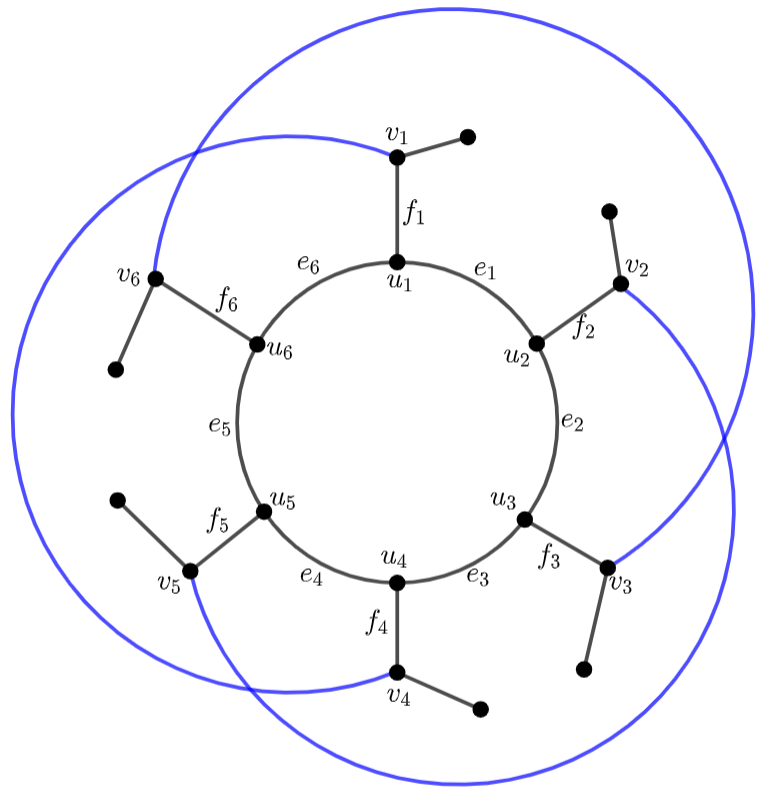}
        \caption{}
        \label{fig
:sub1}
    \end{subfigure}
    \hfill
    \begin{subfigure}[b]{0.3
\textwidth}
        \includegraphics[width=\linewidth]{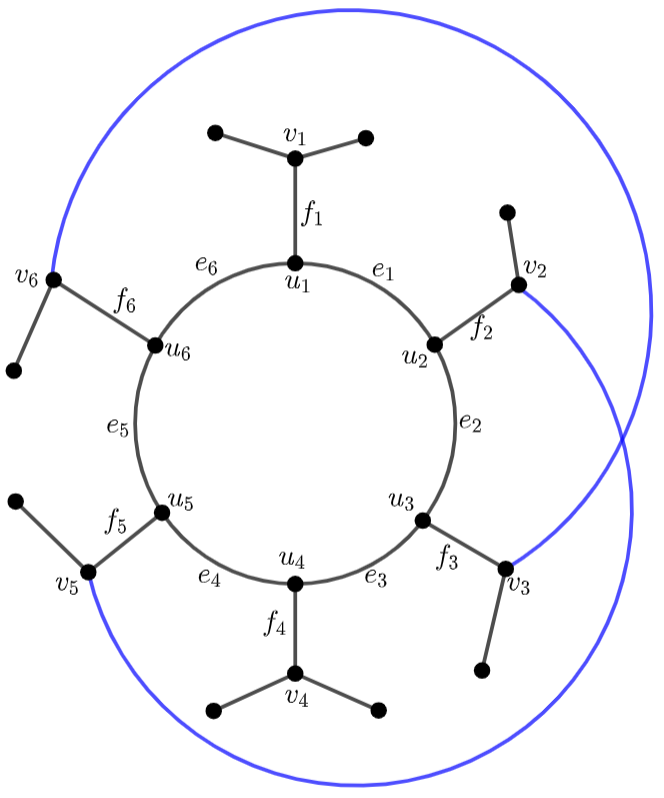}
        \caption{}
        \label{fig
:sub2}
    \end{subfigure}
    \hfill
    \begin{subfigure}[b]{0.3
\textwidth}
        \includegraphics[width=\linewidth]{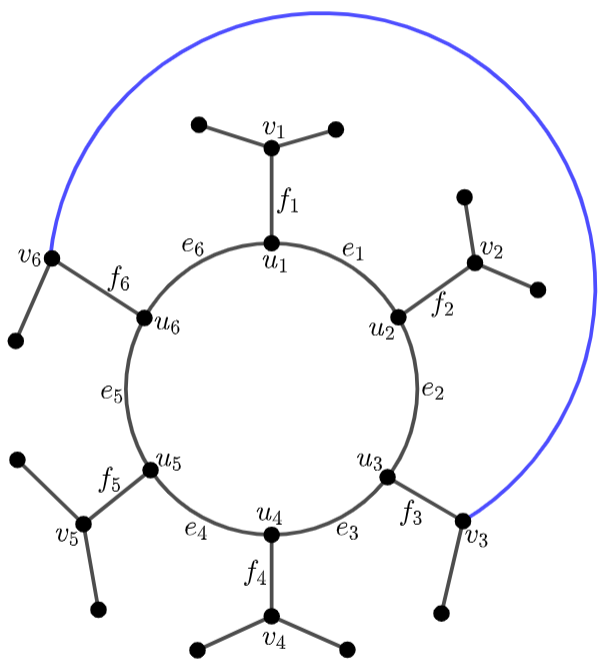}
        \caption{}
        \label{fig
:sub3}
    \end{subfigure}
    \hfill
    \caption{} 

\end{figure}

\begin{Proof}
Suppose that $G$ contains an $6$-cycle $C_6$. Note that $C_6^+$ is an induced subgraph of $G$, and we have three possible local structures based on whether $d_G(f_i,f_j)$ is at most 2 for $i\ne j$.  (see~Figure 1 (1)-(3)). Then for $i\in[6]$, $ |F_\varphi(e_i)| \leq 4$, and $|F_\varphi(f_i)| \leq 5$. Thus $ s(e_i)\geq 6$ and $s(f_i)\geq 5$. For $ i\in[3]$, let $P_i$ be the polynomial $P_H$ of the case in Fig. 1 (i). First we consider $P_1$,
$$P_1=P\left( x_1,x_2,\cdots,x_6,y_1,y_2,\cdots,y_6 \right) =\frac{\prod\limits_{i=1}^{i=3}\left( y_i-y_{i+3} \right) \prod\limits_{1\le i<j\le 6}\left( x_i-x_j \right) \left( x_i-y_j \right) \left( y_i-y_{i+1} \right)}{\prod\limits_{i=1}^{i=3}\left( x_i-x_{i+3} \right)  \prod\limits_{i=1}^{i=6}\left( x_i-y_{i+3} \right) \left( x_i-y_{i+4} \right)}.
$$
Obviously, $\deg(P_1)=45$. By Corollary~\ref{CN1}, we only need to construct suitable monomials $J$ with $\deg(J)=\deg(P_1)$ such that $\eta_J(P_1)\ne 0$. Choose $J=x_1^4 x_2^5 x_3^5 x_4^5 x_5^5 x_6^5 y_1^3 y_2^2 y_3^2 y_4^3 y_5^3 y_6^3$. With the help of MATHEMATICAL, it shows that $ \eta_J(P_1)= -5$. Note that the polynomial $ P_2 $ and $P_3$ are both sub-polynomials of  $P_1$, then  $P_1(\phi) \neq 0$, which implies that $P_2 (\phi)\neq 0$ and $P_3 (\phi)\neq 0$  as well. Thus, a contradiction. The result follows.
\end{Proof}

\begin{Claim}
$G$ has no $n$-cycles, where $n\ge 7$.
\end{Claim}

\begin{figure}[htbp]
\centering
\includegraphics[scale=0.5]{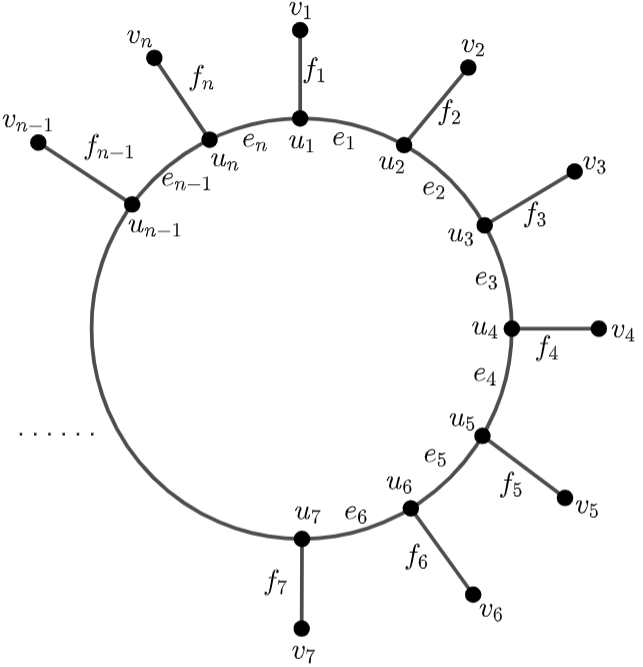}
\caption{}
\label{figure}
\end{figure}

\begin{Proof}
Suppose that $G$ contains an $n$-cycle $C_n$, where $n\ge7$ (see~Figure 2). Since $ |F_\varphi(e_i)| \leq 4$ and $|F_\varphi(f_i)| \leq 6$,  then $ s(e_i) \geq 6$ and $s(f_i) \geq 4$, where $i\in[n]$. For convenience,
let $P = P(x_1, x_2, \cdots, x_n, y_1,  y_2, \cdots, y_n) $, and $x_0=x_n$.
Then
$P=\prod\limits_{i=0}^{i=n}P(i)$, where
\begin{eqnarray*}P(i)&=&
\begin{cases}
\left(\prod\limits_{l=i+1}^{i+2}(x_i-x_l)(x_i-y_l)( y_{i+1}-x_l )\right)(y_{i+1}-y_{i+2}) ~& \mbox{if $0\le i\le n-2$},\\
(x_{n-1}-x_n)(x_{n-1}-y_n)(y_n-x_n)~& \mbox{if $i=n-1$},\\
(x_1-x_{n-1})( x_1-y_n)( y_1-x_{n-1})(y_1-y_n) ~&\mbox{if $i=n$}.
\end{cases}
\end{eqnarray*}
Choose $J=J(0)\prod\limits_{i=5}^{i=n-1}J(i)$, where
\begin{eqnarray*}J(i)&=&
\begin{cases}
x_1^5 x_2^5 x_3^4 x_4^5 y_1^3 y_2^3 y_3^3 y_4^3 y_5^3 ~& \mbox{if $i=0$},\\
x_{i}^{4}y_{i+1}^{3}~& \mbox{if $5\le i\le n-2$},\\
x_{n-1}^{4}y_{n}^{2}x_{n}^{2} ~&\mbox{if $i=n-1$}.
\end{cases}
\end{eqnarray*}
If $n=7$, we can obtain $\eta_J[P]=1$ directly by MATHEMATICAL.

Suppose that $n\ge8$.
With the help of MATHEMATICAL, we will determine $\eta_J[P]$ by some procedures in the following.
$$\eta _{J(0)}\left[ \left(\prod\limits_{0}^{4}P(i)\right)P(n) \right] =(-1)^4x_{n-1}^{2}( x_6+y_6)y_{n}^{2}.$$
For $5\le k\le n-2$,
$$\eta _{J(k)}\left[P(k)\cdot (-1)^{k-1}(x_{k+1}+y_{k+1}) \right]=(-1)^k( x_{k+2}+y_{k+2}).$$
For convenience, for $5\le k\le n-2$, let
$$
\widetilde{P}_{k}=x_{n-1}^{2}y_{n}^{2}\prod\limits_{i=k}^{i=n-1}P(i).
$$
Note that $\widetilde{P}_{k}=\widetilde{P}_{k+1}\cdot P(k)$, $\widetilde{P}_5$ doesn't contain any variable in $J_0$, and $\widetilde{P}_k$ doesn't contain any variable in $J_0\cup_{i=5}^{i=k-1} J_{i}$ for $6\le k\le n-2$. Thus for $5\le k\le n-2$, we have
\begin{eqnarray*}
\eta_{J}[P]
&=&\eta_{J(0)\prod\limits_{i=5}^{i=n-1}J(i)}\left[\prod\limits_{i=0}^{i=n}P(i)\right]\\
&=&\eta_{\prod\limits_{i=5}^{i=n-1}J(i)}\left[[\prod\limits_{i=5}^{i=n-1}P(i)]\cdot \eta_{J_0}[P(n) \prod\limits_{i=0}^{i=4}P(i)]\right]\\
&=&\eta_{\prod\limits_{i=5}^{i=n-1}J(i)}\left[\widetilde{P}_5\cdot(-1)^4(x_6+y_6)\right]\\
&=&\eta_{\prod\limits_{i=k}^{i=n-1}J(i)}\left[\widetilde{P}_k\cdot(-1)^{k-1}(x_{k+1}+y_{k+1})\right]\\
&=&\eta_{\prod\limits_{i=k+1}^{i=n-1}J(i)}\left[\widetilde{P}_{k+1}\cdot \eta_{J_{k}}[P(k)\cdot(-1)^{k-1}(x_{k+1}+y_{k+1})]\right]\\
&=&\eta_{\prod\limits_{i=k+1}^{i=n-1}J(i)}\left[\widetilde{P}_{k+1}\cdot (-1)^{k}(x_{k+2}+y_{k+2})\right]\\
&=&\eta_{J(n-1)}\left[\widetilde{P}_{n-1}\cdot (-1)^{n-2}(x_n+y_n)\right]\\
&=&\eta_{x_{n-1}^{4}y_{n}^{2}x_{n}^{2}}\left[(-1)^{n-2}x_{n-1}^{2} y_{n}^{2}(x_{n-1}-x_n)(x_{n-1}-y_n)(y_n-x_n)(x_n+y_n)\right]\\
&=&(-1)^{n-1}.
\end{eqnarray*}
By Corollary~\ref{CN1}, it is known that when $n \geq 7 $, all $n$-cycles can be strongly $10$-edge colored, a contradiction. Thus the result follows.
\end{Proof}
From the above, we find that a minimal counterexample $G$ does not exist, a contradiction. Thus, we have completed the proof of Theorem \ref{key}.

      \end{document}